\documentclass[12pt]{amsart}

\newcommand{\res}{{\mtr{Res}}}

\newcommand{\barci}{\mtr C_{i}}

\newcommand\numeq[1]%
  {\stackrel{\scriptscriptstyle(\mkern-1.5mu#1\mkern-1.5mu)}{=}}

\newcounter{relctr} 
\everydisplay\expandafter{\the\everydisplay\setcounter{relctr}{0}} 

\newcommand{\barcl}{\mtr C_l}
\newcommand{\barcj}{\mtr C_j}
\newcommand{\barck}{\mtr C_k}

\newcommand{\jo}{{j_1}}
\newcommand{\jtw}{{j_2}}
\newcommand{\barcjo}{\mtr C_{j_1}}
\newcommand{\barcjtw}{\mtr C_{j_2}} 

\usepackage{tikz-cd}
\usepackage{latexsym,amssymb,amsmath}
\usepackage{hyperref}
\usepackage{comment}
\usepackage{mathrsfs}
\usepackage{amsmath,amscd}
\usepackage[enableskew]{youngtab}
\usepackage[all]{xy}
\usepackage{empheq}

\usepackage{amsmath}

\usepackage{empheq}
\usepackage{xcolor}
\definecolor{lightgreen}{HTML}{90EE90}

\newcommand{\cocc}{\co(\cc)}

\newcommand{\sgi}{\sigma(i)}

\newcommand{\mul}{\mu_{\barzm}^l}

 \newcommand{\cfzcc}{{\cf(\czcc)}}
 
 \newcommand{\vsk}{\vskip 0.15cm \noindent}

\newcommand{\ccj}{\cc_j}

 \newcommand{\fpcc}{\fp(\cc)}

\newcommand{\sent}{\mapsto}

\newcommand{\mui}{\mu_i}\newcommand{\muj}{\mu_j}

\numberwithin{equation}{section}

\newcommand\C{\mathcal{C}}
\DeclareMathOperator{\id}{id}

\DeclareMathOperator{\ev}{ev}

\excludecomment{verlong}
\includecomment{vershort}
\excludecomment{noncompile}

\usepackage{tikz}
\usetikzlibrary{matrix}
\usepackage{combelow}

\usepackage{amsmath,amsthm,amsfonts,amssymb}
\usepackage{amsmath, amsthm, amssymb, amscd, amsfonts}
\usepackage[all]{xy}
\usepackage{amssymb, amsthm, amsmath}
\usepackage{amsfonts}
\usepackage{color}

\newcommand{\ra}{\rightarrow}

\newcommand{\Z}{\mathbb Z}

\newcommand{\ot}{\otimes}
\newcommand{\co}{\mathcal O}
\newcommand{\xra}{\xrightarrow}
\newcommand{\mtc}{\mathcal}

\newcommand{\lam}{\lambda}

\newcommand{\Lam}{\Lambda}

\newcommand{\al}{\alpha}
\newcommand{\eps}{\epsilon}

\newcommand{\ul}{\underline}

\newcommand{\lh}{\leftharpoonup}
\numberwithin{equation}{section}
\newtheorem{lem}[equation]{Lemma}

\theoremstyle{plain}

\newtheorem{thm}[equation]{Theorem}
\newtheorem{prop}[equation]{Proposition}
\newtheorem{defn}[equation]{Definition}
\newtheorem{cor}[equation]{Corollary}
\newtheorem{rem}[equation]{Remark}
\newcommand{\ch}{\chi}
\newcommand{\mtr}{\mathrm}

\numberwithin{equation}{section}
\newcommand{\ncm}{\newcommand}
\ncm{\np}{\newpage}
\ncm{\ebl}{\end{thebibliography}}
\ncm{\bbl}{\begin{thebibliography}}
\ncm{\chd}{_{ _{\ch}}}
\ncm{\ald}{_{ _{\al}}}
\newcommand{\blam}{\Lam}
\ncm{\cP}{\mathcal{P}}
\ncm{\ei}{e_i}
\ncm{\eij}{e_{i,\;j}}
\ncm{\bt}{\begin{thm}}
\ncm{\bdef}{\begin{defn}}
\ncm{\edf}{\end{defn}}
\ncm{\et}{\end{thm}}
\ncm{\bc}{\begin{cor}}
\ncm{\bl}{\begin{lem}}
\ncm{\el}{\end{lem}}
\ncm{\bpf}{\begin{proof}}
\ncm{\epf}{\end{proof}}
\ncm{\ec}{\end{cor}}
\ncm{\ord}{\mtr{ord}}
\ncm{\er}{\end{rem}}
\ncm{\br}{\begin{rem}}
\ncm{\bn}{\begin}

\ncm{\bp}{\begin{prop}}
\ncm{\ep}{\end{prop}}
\ncm{\bd}{
\begin{document}}
\ncm{\ed}{\end{document}}
\ncm{\beq}{\begin{equation}}
\ncm{\beqn}{\begin{equation*}}
\ncm{\eeq}{\end{equation}}
\ncm{\eeqn}{\end{equation*}}
\ncm{\bea}{\begin{eqnarray}}
\ncm{\eea}{\end{eqnarray}}
\ncm{\beanon}{\begin{eqnarray*}}
\ncm{\eeanon}{\end{eqnarray*}}\ncm{\ek}{\eps|_K}\ncm{\diez}{\#}
\ncm{\bwt}{\bowtie}
\ncm{\cC}{\mtc{C}}\ncm{\cc}{\mtc{C}}
\ncm{\cX}{\mtc{X}}
\ncm{\wt}{\widetilde}
\ncm{\sg}{\sigma}
\ncm{\Rep}{\mathrm{Rep}}
\DeclareMathOperator{\Irr}{Irr}
\ncm{\X}{\mathcal{X}}
\ncm{\cA}{\mathcal{A}}
\ncm{\HKer}{\mtr{HKer}}
\ncm{\LKER}{\mtr{LKer}}
\ncm{\aad}{\mtr{ad}}
\newcommand{\mbf}{\mathbb F}
\ncm{\Dr}{\mtr{D}}
\ncm{\cD}{{\mathcal{D}}}\ncm{\cd}{{\mathcal{D}}}\ncm{\ce}{{\mathcal{E}}}
\ncm{\G}{\mathcal{G}}
\ncm{\Dc}{\mtc{D}}
\ncm{\E}{\mtc{E}}
\ncm{\fp}{\mtr{FPdim}}
\ncm{\Vc}{\mtr{Vec}}
\ncm{\cK}{\mtc{K}}
\ncm{\cM}{\mtc{M}}
\ncm{\cE}{\mtc{E}}
\ncm{\cS}{\mtc{S}}

\newcommand{{\ipr}}{i'}

\DeclareMathOperator{\End}{End}
\ncm{\cop}{\mtr{cop}}
\ncm{\op}{\mtr{op}}
\ncm{\chr}{character }\ncm{\ck}{\mtc{K}}
\ncm{\bw}{\bwt}
\ncm{\hker}{\mtr{HKer}}
\ncm{\bx}{\boxtimes}
\ncm{\blue}{\textcolor[rgb]{.00, .00, 1.00}}
\ncm{\red}{\textcolor[rgb]{1.00, .00, .00}}
\ncm{\green}{\textcolor[rgb]{.50, 0.20, .90}}
\ncm{\bne}{\begin{enumerate}}
\ncm{\ene}{\end{enumerate}}
\ncm{\lker}{\mtr{LKer}}
\ncm{\md}{\medbreak}
\ncm{\rep}{\Rep}\ncm{\ind}{\mtr{ind}}
\ncm{\mdn}{\md\noindent}
\ncm{\dd}{$}
\ncm{\up}{^}
\newcommand{\tcs}{\text}
\newcommand{\mbb}{\mathbb B}
\newcommand{\vs}{\mathbb V}
\newcommand{\sth}{suppose that\;}
\newcommand\rad{\operatorname{rad}}
\newcommand{\itm}{\item}
\newcommand{\dbd}{$$}
\newcommand{\mol}{\mtr{mod}}
 \newcommand{\ro}{\rho}
\newcommand{\irr}{\mathrm{Irr}}
\newcommand{\mbc}{\mathbb C}
\newcommand{\mbs}{\mathbb S}
\newcommand{\mbz}{\mathbb Z}
\newcommand{\ct}{\mtc T}
\newcommand{\sm}{\setminus}
\newcommand{\epl}{^{+}}
\newcommand{\sbsq}{\subseteq}
\newcommand{\sbs}{\subset}
\newcommand{\cco}{\mtr{co}}
\newcommand{\cz}{\mathcal{Z}}
\newcommand{\dual}{^{*}}
\newcommand{\Gm}{\Gamma}
\ncm{\cY}{\mtc{Y}}
\newcommand\ZZ{{\mathbb Z}} 
\newcommand{\bab}{\color{DarkOrchid}{}}
\newcommand{\eab}{\normalcolor{}}
\newcommand{\subs}{\subsection}
\newcommand{\cv}{\mtc{V}}
  \newcommand{\grn}{\green}
\newcommand{\dt}{\delta}

\newcommand{\ccf}{\mathrm{ {CF}(\cc)}}
\newcommand{\cce}{\mathrm{ {CE}(\cc)}}
\newcommand{\cecc}{\mathrm{ {CE}(\cc)}}
\newcommand{\cecd}{\mathrm{ {CE}(\cd)}}
\newcommand{\kk}{\Bbbk}
\newcommand{\otL}{\ot_{L}}
\newcommand{\otl}{\ot_{L}}
\newcommand{\unpsi}{1_{\psi}}
\newcommand{\epsi}{e_{\psi}}
\newcommand{\ephi}{e_{\phi}}
\newcommand{\ech}{e_{\ch}}
\newcommand{\nleftcid}{\text{left normal  coideal subalgebra}}
\newcommand{\dimL}{\dim_{\kk}L}
\newcommand{\cl}{\mtc L}
\newcommand{\mj}{\mtc J}
\newcommand{\tl}{\tilde L}
\newcommand{\tL}{\tilde L}
\newcommand{\tpsi}{\tilde(\psi)}
\newcommand{\tmx}{\tilde{\mtc X}}
\newcommand{\zlh}{\mathrm{ZL}}
\newcommand{\ba}{\mathrm A}
\newcommand{\bv}{\mathrm V}
\newcommand{\zhopf}{\mtc{Z}_{\mtr{Hopf}}}
\newcommand{\lstar}{L^{*}}
\newcommand{\ldstar}{L^{**}}
\newcommand{\mstar}{M^{*}}
\newcommand{\mdstar}{M^{**}}
\newcommand{\lkera}{\lker_{A}}
\newcommand{\mdprime}{M''}
\newcommand{\ldprime}{L''}
\newcommand{\cm}{\mtc M}
\newcommand{\ccm}{\mathcal M}
\newcommand{\cn}{\mathcal N}
\newcommand{\ccn}{\mathcal N}
\newcommand{\rx}{\mtr{Rex}}
\newcommand{\cca}{\ca}
\newcommand{\ih}{\underline{\mtr{Hom}}}
\newcommand{\cih}{\underline{\mtr{coHom}}}
\newcommand{\hm}{\mtr{ {Hom}}}
\newcommand{\cov}{\mtr{coev}}
\newcommand{\rora}{\rho^{\mtr{ra}}}
\newcommand{\rola}{\rho^{\mtr{la}}}
\newcommand{\cx}{\mtc X}
 \newcommand{\cZ}{\cz}
 \newcommand{\ca}{\cA}
 \newcommand{\stat}{\noindent}
 \newcommand{\bfa}{{\bf A}}
 \newcommand{\unu}{\mathbf{1}}
 \newcommand{\barzu}{{\bar {  Z}(\unu)}}
 
\newcommand{\idx}{\id_X}
\newcommand{\lprime}{L'}
\newcommand{\mprime}{M'}
\newcommand{\nat}{ \mtr{{  Nat}}}
\newcommand{\ft}{\mtc F_\lam}
\newcommand{\rhau}{\rightharpoonup}
\newcommand{\lhau}{\leftharpoonup}
\newcommand{\cf}{\mathrm{ {CF}}}

\newcommand{\cfc}{\mathrm{{CF}}(\cc)}
\newcommand{\csu}{\overline{\mathfrak{  C}}}
\newcommand{\cfcc}{\mathrm{ {CF}}(\cc)}
\newcommand{\cfcd}{\mathrm{CF}(\cd)}
\newcommand{\cfd}{\mathrm{CF}(\cd)}
\newcommand{\czcc}{{\cz(\cc)}}
\newcommand{\czcd}{{\cz(\cd)}}
\newcommand{\czt}{{\cz(\cz(\cc))}}
\newcommand{\enx}{\mtr{  End}}
\newcommand{\runu}{R(\unu)}

\newcommand{\bdfn}{\bn{defn}}
\newcommand{\edfn}{\end{defn}}
\newcommand{\deltax}{\delta_X}
\newcommand{\deltav}{\delta_V}
\newcommand{\repcca}{\rep_\cc(A)}
\newcommand{\xotay}{X \ot_A Y}
\newcommand{\xoty}{X \ot Y}
\newcommand{\votw}{V \ot W}
\newcommand{\votaw}{V \ot_A W}
\newcommand{\dimax}{\dim_AX}
\newcommand{\dimccx}{\dim_\cc(X)}
\newcommand{\dimcca}{\dim_\cc(A)}
\newcommand{\dimccv}{\dim_\cc(V)}
\newcommand{\dima}{\dim_A}
\newcommand{\biga}{A}
\newcommand{\comp}{\mathbb C}
\newcommand{\tehtaa}{\theta_A}
\newcommand{\tetaa}{\theta_A}
\newcommand{\ida}{\id_A}
\newcommand{\hma}{\hm_A}
\newcommand{\hmcc}{\hm_\cc}
\newcommand{\fv}{F(V)}
\newcommand{\fw}{F(W)}
\newcommand{\ota}{\ot_A}
\newcommand{\repza}{\rep_\cc^0(A)}
\newcommand{\epsa}{\eps_A}
\newcommand{\bndefn}{\bn{defn}}
\newcommand{\edefn}{\end{defn}}
\newcommand{\bdefn}{\bn{defn}}

\newcommand{\vld}{V^{*}}
\newcommand{\vldd}{V^{**}}
\newcommand{\xld}{X^{*}}
\newcommand{\xldd}{X^{**}}
\newcommand{\yld}{Y^{*}}
\newcommand{\yldd}{Y^{**}}
\newcommand{\aldu}{A^{*}}
\newcommand{\aldd}{A^{**}}

\newcommand{\ia}{\mtr{i}_A}
\newcommand{\aota}{A\ot A}

\newcommand{\idv}{\id_V}

\newcommand{\ld}{^*}
\newcommand{\repg}{\rep(G)}

\newcommand{\thetav}{\theta_V}

\newcommand{\tta}{\theta_A}

\newcommand{\muv}{\mu_V}
\newcommand{\muw}{\mu_W}

\newcommand{\dimcc}{\dim(\cc)}
\newcommand{\chii}{\chi_i}
\newcommand{\chistar}{\ch_{i^*}}
\newcommand{\chj}{\ch_j}
\newcommand{\chm}{\ch_m}
\newcommand{\chn}{\ch_n}
\newcommand{\dimvi}{\dim(V_i)}
\newcommand{\mtcd}{Q}
\newcommand{\mtca}{\mtc A}
\newcommand{\lamcd}{\lam_\cd}
\newcommand{\fpdimcd}{\fp(\cd)}
\newcommand{\laml}{\lam_L}
\newcommand{\apm}{A//M}
\newcommand{\apl}{A//L}
\newcommand{\repapm}{\rep(\apm)}
\newcommand{\repapl}{\rep(\apl)}
\newcommand{\dimvj}{\dim(V_j)}
\newcommand{\dvi}{\dim(V_i)}
\newcommand{\dvj}{\dim(V_j)}
\newcommand{\sumjtom}{\sum_{j=0}^m}
\newcommand{\sumitom}{\sum_{i=0}^m}
\newcommand{\sij}{s_{ij}}
\newcommand{\sji}{s_{ji}}
\newcommand{\dxj}{d_j}
\newcommand{\dxi}{d_i}
\newcommand{\dimka}{\dim_{\kk}(A)}
\newcommand{\dimk}{\dim_{\kk}}
\newcommand{\blaml}{\blam_L}
\newcommand{\sumjtor}{\sum_{j=0}^r}
\newcommand{\dimkl}{\dim_{\kk}(L)}
\newcommand{\mtcjl}{\mtc J_L}
\newcommand{\vota}{ V\ot A}
\newcommand{\vi}{V_i}
\newcommand{\vj}{V_j}
\newcommand{\dimcd}{\dim(\cd)}

\newcommand{\alij}{\al_{ij}}
\newcommand{\alji}{\al_{ji}}
\newcommand{\rcc}{r_\cc}
\newcommand{\rcd}{r_\cd}
\newcommand{\clsx}{[X]}
\newcommand{\clsy}{[Y]}
\newcommand{\clsz}{[Z]}
\newcommand{\rcdp}{r_{\cd'}}
\newcommand{\sumjtorp}{\sum_{j=0}^{r'}}
\newcommand{\aljm}{\al_{jm}}
\newcommand{\aljn}{\al_{jn}}
\newcommand{\sjm}{s_{jm}}
\newcommand{\smj}{s_{mj}}
\newcommand{\snj}{s_{nj}}

\newcommand{\betaij}{\beta_{ij}}
\newcommand{\betaji}{\beta_{ji}}

 \newcommand{\ip}{i'}
\newcommand{\sumjtoprp}{\sum_{j=0}^{r'}}
\newcommand{\sumjtopr}{\sum_{j=0}^{r}}
 \newcommand{\teh}{\tilde{h}}
\newcommand{\cdp}{\cd'}
\newcommand{\xphii}{X_{\phi(i)}}
\newcommand{\inv}{^{-1}}

\newcommand{\fq}{f_Q}
\newcommand{\tr}{\mtr{tr}}
\newcommand{\rtwone}{R_{21}R}

\newcommand{\ccad}{{\cc_{\mtr{ad}}}}
\newcommand{\ccpt}{{\cc_{\mtr{pt}}}}
\newcommand{\qtr}{quasi-triangular\;}
\newcommand{\trq}{\tr_q}

\newcommand{\repal}{\mtr{Rep}(A//L)}
\newcommand{\lkeravi}{\lker_A(V_i)}
\newcommand{\lkeravj}{\lker_A(V_j)}
\newcommand{\cross}[1][1pt]{\ooalign{%
 \rule[1ex]{1ex}{#1}\cr
 \hss\rule{#1}{.7em}\hss\cr}}
\newcommand{\blml}{\blam_L} 
\newcommand{\phir}{\phi_R}
\newcommand{\kda}{{  \Phi(A)}}

\newcommand{\mtcil}{\mtc{I}_L}

\newcommand{\un}{\unu}
\newcommand{\tfl}{\mtc{T}}
\newcommand{\barzm}{\barz(M)}
\newcommand{\barzn}{\barz(N)}
\newcommand{\ccr}{\mtc R^{\cc}}
\newcommand{\ulc}{\ul{\cc}}

\newcommand{\pimx}{\pi_{M;\;X}}
\newcommand{\pinx}{\pi_{N;\;X}}
\newcommand{\acc}{{\mathrm A_\cc}}
\newcommand{\epsu}{\eps_\unu}

\newcommand{\ob}{\mtr{Obj}}
\newcommand{\obc}{\mtr{Obj(\cc)}}
\newcommand{\ccop}{\cc^{\mtr{op}}}
\newcommand{\mtf}{\mtc F_\lam}
\newcommand{\mtfi}{\mtc F^{-1}_\lam}
\newcommand{\elcd}{\ell_\cd}
\newcommand{\mcid}{\mtc I_\cd}
\newcommand{\mcidp}{\mtc I_{\cd'}}
\newcommand{\wtildelcd}{\widetilde{\elcd}}
\newcommand{\wtildelcdp}{\widetilde{\ell_{\cd'}}}
\newcommand{\cpt}{\cc_{\mtr{pt}}}
\newcommand{\barzr}{\barz_\cd}
\newcommand{\barzv}{\barz(V)}
\newcommand{\acd}{\mathrm A_\cd}
\newcommand{\czrcd}{\cz_\cc(\cd)}
\newcommand{\sml}{\Small}
\newcommand{\bs}{\blue{\Small }}
\newcommand{\yd}{Yetter-Drinfeld}

\newcommand{\sumitor}{\sum_{i=0}^r}
\newcommand{\cdop}{\cd^{\mtr{op}}}
\newcommand{\ccrev}{\cc^{\mtr{rev}}}
\newcommand{\barz}{{\bar{\mathrm Z}}}
\newcommand{\etl}{etale\;}
\newcommand{\czca}{\cz(\ca)}

\newcommand{\tq}{\tilde Q}
\newcommand{\tqu}{\tilde Q^1}
\newcommand{\tqt}{\tilde Q^2}
\newcommand{\phip}{\phi'}
\newcommand{\rphip}{\;_R\phi'}
\newcommand{\phipr}{\phi'_R}
\newcommand{\tlr}{\tilde{R}}
\newcommand{\tlrm}{\tilde{r}}

\newcommand{\tetx}{\text}
\newcommand{\widehta}{\widehat}
\newcommand{\wdhat}{\widehat}
\newcommand{\wht}{\widehat}
\newcommand{\cofa}{{\mathbb C[A]}}
\newcommand{\wdt}{\widehat}
\newcommand{\dl}{{^*}}
\newcommand{\comx}{\mathbb C}
\newcommand{\ovl}{\overline}

\newcommand{\mujo}{\mu_\jo}
\newcommand{\mujtw}{\mu_\jtw}
\newcommand{\adz}{a{\dl}}
\newcommand{\bdz}{b{\dl}}

\newcommand{\spr}{S^\perp}
\newcommand{\cofs}{\comp [S]}
\newcommand{\spz}{S^{\perp_z}}

\newcommand{\omz}{\omega_z}
\newcommand{\zg}{\mathrm{Z}(S)}
\newcommand{\aling}{{\al \in g}}

\newcommand{\blkg}{\mtr{Bl}(g)}
\newcommand{\clsg}{\mtr{Cl}(g)}
\newcommand{\mtadinv}{(\mtc G^{{\wdht A}})^{{-1}}}
\newcommand{\muk}{\mu_{k}}
\newcommand{\mta}{\mtc F^{A}}
\newcommand{\cofad}{\comp[\wdht A]}
\newcommand{\wtau}{\wdht{\tau}}
\newcommand{\mtainv}{{(\mta)}^{-1}}
\newcommand{\wdht}{\widehat}
\newcommand{\augm}{\mtr{aug}}
\newcommand{\mua}{\wdht {\wdht a}}
\newcommand{\aps}{A//S}
\newcommand{\cssa}{\cc(S, A)}
\newcommand{\aug}{\mtr{aug}}
\newcommand{\rss}{{\big|_S}}
\newcommand{\gprp}{g^\perp}
\newcommand{\alins}{{s \in S}}

\newcommand{\sz}{\widehta{s}}
\newcommand{\wmu}{\widehta{\mu}}
\newcommand{\wmui}{\widehta{\mu}_i}
\newcommand{\wmuj}{\widehta{\mu}_j}
\newcommand{\wch}{\widehta{\ch}}

\newcommand{\clsc}{\mtr{Cl}(\cc)}

\newcommand{\wf}{\widehat{F}}
\newcommand{\whta}{\widehat A}
\newcommand{\wtf}{\widehat{F}}
\newcommand{\sgj}{{\sg(j)}}
\newcommand{\wzd}{\widehat{d}}
\newcommand{\wpm}{\widehat{P}}
\newcommand{\wps}{\widehat{p}}

\newcommand{\oic}{\overline{\irr(\cc)}}
\newcommand{\clc}{\overline{\mtr{Cls}(\cc)}}

\newcommand{\wh}{\widehat{h}}
\newcommand{\rgo}{{\mathbb R_{\geq 0}}}
\newcommand{\sumini}{\sum_{i\in I}}
\newcommand{\cofb}{\comp[B]}

\newcommand{\sgwhf}{\sg_{_{\widehat F}}}

\newcommand{\we}{\widehta{E}}
\newcommand{\sumktom}{\sum_{k=0}^m}
\newcommand{\gal}{\mtr{Gal}}
\newcommand{\galkq}{\gal(\mathbb K/\mathbb Q)}
\newcommand{\sgf}{\sg_{_F}}
\newcommand{\sggi}{{\sg(i)}}
\newcommand{\sge}{\sg_{_E}}
\newcommand{\unue}{{1_{\cecc}}}
\newcommand{\hsgj}{\widehat{\sg}(j)}
\newcommand{\whsgi}{\widehta{\sg}(i)}

\newcommand{\wpp}{\widehat{p}}
\newcommand{\dimccj}{\dim(\cc^j)}
\newcommand{\tauj}{{\tau(j)}}
\newcommand{\dimcctauj}{\dim(\cc^\tauj)}
\newcommand{\etas}{{\eta(s)}}
\newcommand{\mcc}{m_\cc}

\newcommand{\wal}{\widehta{\al}}
\newcommand{\wj}{\widehat{\mtc J}}
\newcommand{\galc}{\mtr{Gal}_{\cc}}
\newcommand{\galz}{\mtr{Gal}_{\czcc}}
\newcommand{\wjr}{\widehat{J}_{R}}

\newcommand{\dimcck}{\dim(\cc^k)}
\bd
\title[Fusion categories]
{Structure constants for premodular categories}
\author{Sebastian Burciu}
\address{Inst.\ of Math.\ ``Simion Stoilow" of the Romanian Academy P.O. Box 1-764, RO-014700, Bucharest, Romania}
\email{sebastian.burciu@imar.ro}
\date{\today}
\maketitle
\section*{Abstract}
In this paper we study conjugacy classes for  pivotal  fusion categories. In particular we prove a Burnside type formula for the structure constants concerning the product of two conjugacy class sums of a such fusion category.  For a braided weakly integral    fusion category $\cc$ we show that these structure constants multiplied by $\dimcc$ are non-negative integers, extending some results obtained by Zhou and Zhu (see \cite{zz}) in the settings of  semisimple quasitriangular Hopf algebras.  
\section{Introduction}
Conjugacy classes for finite groups are a very important tool in the study of representations of finite groups. Their associated class sums form a basis for the center of the group algebra and the structure constants obtained by multiplying two such class sums are integers. A precise formula for these structure constants, in terms of characters, was given by Burnside,  see e.g. \cite[pp. 45]{ism}.%

In the spirit of Zhu's work, \cite{zind}, Cohen and Westreich in \cite{CW6} have defined a notion of  a conjugacy class of a semisimple Hopf algebra. Using this new concept they were able to extend some results from finite group representations to semisimple Hopf algebra representations. 

More recently, Shimizu  introduced  in \cite{scalg},  the notion of conjugacy classes for fusion categories, extending the previous notion of  Cohen and Westreich  mentioned above.  In the same paper Shimizu associated to each conjugacy class a central element called also conjugacy class sum.  These class sums play the role of the sum of group elements of a conjugacy class in group theory.

In \cite{scalg} Shimizu asked what results from \cite{CW6, CW4} can be extended from semisimple Hopf algebras to fusion categories.

In this paper we give an analogue of Burnside formula for the structure constants of any pivotal fusion category extending the results obtained in \cite{CW6} for semisimple Hopf algebras:
\bt\label{burnside} 
Let $\cc$ be a pivotal fusion category with a commutative Grothendieck ring. With the above notations one has 
\beq\label{bnf}
c^k_{ij}=\sum_{s=0}^m \frac{\ch_s( \barci)\ch_s(\barcj)\ch_{s^*}(\barck)}{\dimcc \dim(\cc^k)d_s}.
\eeq
\et
To prove the above theorem, we use  the theory of probability groups, introduced by Harrison in 1979 in \cite{hdk}, as a tool to study the character ring of a pivotal fusion category. Extending the results from \cite{zz}, for a pivotal fusion category $\cc$, we prove that the algebra generated by the dual of the probability group $\cfcc$ is isomorphic to the center $\cecc$ as introduced by Shimizu in \cite{scalg}, see also Section \ref{fcgg}.

Cohen and  Westreich \cite[Theorem 2.6]{CW6} proved that the product of two class
sums of H is an integral combination up to a factor of $\dim(H)^{-2}$ of the class sums of H, i.e.
$$
\barci\barcj =\frac{1}{\dim(H)^2}\big(\sumktom \widetilde{c}^{\;k}_{ij} \barck\big)
$$
where $\{\barci \;| \;0 \leq i\leq  m\}$ is the set of class sums of $H$ and each $
\widetilde{c}^{\;k}_{ij}$ is a non-negative integer. Recently in \cite[Theorem 5.6]{zz} the authors have shown that this factor can be replaced by $\dim(H)$. Next theorem generalizes the result obtained by Zhu and Zhou in \cite{zz} from quasi-triangular Hopf algebras to integral premodular  categories.
\bt\label{main2}
Let $\cc$ be a  premodular  category. Then the structure constants $c^l_{j_1, j_2}$ multiplied by $\dim(\cc)$ are algebraic integers. 

Moreover, if $\cc$ is weakly integral, then the same numbers multiplied by $\dim(\cc)$ are  non-negative integers.
\et
The proof of theorem has as a key step a result of Ostrik from \cite[Theorem 2.13]{O3}, concerning the image by the forgetful functor of the central primitive idempotents of the Drinfeld center $\czcc$ from $\cf(\czcc)$ to $\cfcc$.

We also show that the character ring $\cfcc$ of a modular fusion category $\cc$ is self dual, generalizing the result  obtained by Zhu and Zhou in \cite{zz} for the Drinfeld double $D(H)$ of a semisimple Hopf algebra $H$.

This paper is organized as follows. In Section \ref{pg} we recall the basic facts on probability groups and prove a formula for the (hypergroup) algebra structure of the dual of a probability group. In Section \ref{fcgg} we recall the basics on fusion categories and the describe the dual probability group structure of the character ring as the space of central elements. Theorem \ref{burnside} is proven in Subsection \ref{proof}. Section \ref{braided} contains the proof of Theorem \ref{main2}. 

We work over the ground filed $\comp$ of complex numbers. All algebras and fusion  categories under consideration are $\comp$-linear.
\section{Probability groups} \label{pg}
First we recall the definition of a probability group as in \cite{hdk}.
A {\it probability goup} is a set $A$ together with a function 
$$
p:A^3\ra \mathbb R_{\geq 0}, \;(a,b,c)\sent p(a.b=c)
$$ 
satisfying the following conditions:
\bne
\item
 For all $a, b\in A$, the values $p(a.b=c)=0$ for all but finitely many $c \in A$, and
 \beq\label{c1}
p(a.b=c)\geq 0,\;\sum_{c \in A}p(a.b=c)=1
\eeq
\item For all $a,b, c,d\in A$ one has

\beq\label{c5-assoc}
\sum_{x \in A}p(ab=x)p(xc=d)=\sum_{y \in A}p(ay=d)p(bc=y)
\eeq
\item There is an identity  element $1\in A$ such that for any $a \in A$
\beq\label{c3-unit}
p(1.a=a)=p(a.1=a)=1
\eeq
\item For any $b\in A$  there is a unique $b\dl\in A$ such that
\beq\label{c2}
p(b.b\dl=1)>0
\eeq
\item For any $a,b,c\in A$ one has that 
\beq\label{c3}
p(a.b=c)=p(b\dl .a\dl=c\dl).
\eeq
\item For all $a \in A$ it follows 
\beq\label{c4}
p(a.a\dl=1)=p(a.\dl a=1).
\eeq
\ene
\noindent We also use the notation $p_c(a,b)$ for $p(a.b=c)$. If (A, p) is a probability group, one easily checks that the identity element $1$ is unique, $1\dl=1$ and $(a\dl)^\dl=a$ for all $a \in A$.
\vsk
Given a set $A$ we denote by $\cofa$ the complex vector space with linear basis $a \in A$. Then $(A,p)$ is a probability group if and only if $\cofa$ is a $\kk$-algebra with multiplication given by 
$$
a.b=\sum_{c \in A}p(ab=c)c
$$ 
for all $a,b \in A$.  Note that the associativity of the product as defined above is equivalent to Equation \eqref{c5-assoc}.

Based on the above properties of $p$ it is easy to see that 
\beq\label{ua}
u_A:=\frac{1}{n(A)}\sum_{a\in A} h_aa
\eeq
is an idempotent element of $A$ and $au_A=u_A=u_Aa$ for any $a\in A$.
\newcommand{\aina}{{a\in A}}

We define a bilinear map $m: \cofa\times \cofa\ra \mathbb C$  where  $m(a, x)$ is the coefficient of $a$ in writing $x$ as linear combination of the standard basis $A$,  for any $a \in A$ and $x \in \cofa$. Therefore $m(a, bc)=p_a(b,c)$ for all $a,b,c \in A$.  Extending this by linearity, in general one has that
\beq\label{m}
m_A(\sum_{\aina}\al_aa, \sum_\aina \beta_aa)=\sum_\aina\al_a\beta_a.
\eeq

Recall that a probability group is called {\it abelian} if $p_c(a,b)=p_c(b,a)$ for all $a,b,c \in A$. This equivalent to $\cofa$ being a commutative algebra.

For the rest of this section we let $A$ be a finite abelian probability group  with the cardinal $|A|=m+1$.

Define $\widehat A:=\{\muj:\cofa \ra \comp\}$ the set of all algebra algebra morphisms of $\cofa$. Since $A$ is abelian and $\cofa$ is a semisimple commutative algebra  (see \cite{hdk}) it follows that $\whta$ is a $\comp$-basis of $\cofa^*$. Therefore one has $|A|=|\widehta A|$. Moreover, one can define a multiplication on $\whta$ by declaring
\beq
[\mui\star \muj](a):=\mui(a)\muj(a), \;\text{for all}\;a\in A.
\eeq

Extending linearly $\mui\star \muj$ on the entire $\cofa$ one obtains an algebra structure on $\cofa^*=\comp[\widehta A]$. It follows that there are some non-zero scalars ${\widehat p}_k(i,j)\in \comp$ such that 
\beq\label{hatp}
\mui \star \muj=\sumktom{\wdht p}_k(i,j)\muk.
\eeq 
\subsection{}
Recall \cite{wd} that a {\it generalized hypergroup} is a set $A$ such that $\mathcal A:=\cofa$  is a *-algebra with unit $a_0$ over $\mathbb C$ and $A =\{a_0,a_1\dots,a_n\}$ is a basis of $\mathcal A$ with $\mathcal A ^*=\mathcal A$ for which the structure constants $n_{ij}^k$ defined by $$c_ic_j=\sum_k n_{ij}^k c_k$$ satisfy the  following two conditions: $$c_i^*=c_j \iff n_{ij}^0>0, \;\text{and}\;c_i^*\neq c_j \iff n_{ij}^0 =0.$$

A generalized hypergroup $A$ is called {\it abelian} if $c_ic_j=c_jc_i$ for all $i,j$; it is called {\it real} if $n_{ij}^k\in \mathbb R$ for all $i,j,k$; {\it  positive} if $n_{ij}^k\geq 0$ for all $i,j,k$ and {\it normalized} if $\sum_k n_{ij}^k =1$ for all $i,j$. Therefore a probability group is a real normalized generalized hypergroup with $n^0_{ii*}=n^0_{ii^*}$ for all $i$.


Let $A$ be an abelian probability group. Define the linear functional $\tau:\cofa \ra \comp, \;a\sent \delta_{a,1}$. Note that $\tau(x)=m_A(1,x)$ for all $x \in \cofa$. Since $\tau:\cofa \ra \comp$ is a non-degenerate trace there are some non-zero complex numbers $n_j$  such that $\tau(F_j)=\frac{1}{n_j}$.  
\noindent It follows that
\beq\label{tform}
\tau=\sumjtom\frac{\mu_j}{n_j}
\eeq
 and both $\{a, h_aa\dl\}$ and $\{F_j, \frac{1}{n_j}F_j\}$ are dual bases for $\tau$. Therefore 
\beq\label{clseq}
\sumjtom n_jF_j\ot F_j=\sum_{a\in A} a\ot h_aa\dl
\eeq


Note that $\aug:\cofa\ra \comp,\;a\sent 1$ is a morphism of $\comp$-algebras. Without loss of generality we may suppose that $\mu_0=\aug$. It follows that $\mu_0\star \mu_j=\mu_j\star \mu_0=\mu_j$ for all $j$. Moreover, the central primitive idempotent $F_0$ associated to $\mu_0$ coincides to the element $u_A$ from Equation \eqref{ua}. Moreover, note that $n_0=n(A)$ since
$$
\frac{1}{n_0}=\tau(F_0)=\frac{1}{n(A)}\big(\sum_{a \in a} h_a\tau(a)\big)=\frac{h_1}{n(A)}=\frac{1}{n(A)}.
$$


Let $A$ be an abelian probability group as above and $F_j$ a primitive central idempotent corresponding to a character $\muj:\cofa \ra \comp$.
Applying $\id\ot \mu_j$ to  Equation \eqref{clseq} one obtains:
\beq\label{fj1}
F_j=\frac{1}{n_j}\bigg(\sum_{a \in A}h_a\mu_j(a\dl)a\bigg).
\eeq


Given an algebra homomorphism $\muj:\cofa\ra \comp$  as above, it is easy to see that $\muj^{\dl}:\cofa\ra \comp,\;a\sent \ovl{\muj(a)}$ is also an algebra homomorphism. We let $j\dl$ be the index for which $\mu_j^\dl=\mu_{j\dl}$. With the above notations by \cite[Proposition 2.10]{hdk} it follows that $p_{0}(j, j\dl)>0,\;\tetx{for all}\; 0\leq j\leq m.$
We let $\wh_j:=\frac{1}{p_{0}(j, j\dl)}$ and $n(\widehat A):=\sumjtom \wh_j$.


Following the same \cite[Proposition 2.10]{hdk} for any finite abelian probability group $A$ one has that 
\bne
\item $\wh_j$ is a real positive number.
\item $n(A)=n(\widehat{A})$.
\item 
For any $\jo, \jtw$ the first orthogonality relation is written as:
\beq\label{so1}
\sum_{a \in A}h_a\wh_{\jo}\mu_\jo(a)\mu_\jtw(a\dl)=\delta_{\jo, \jtw}n(A)
\eeq
\item 
For any $a,b \in A$ the second orthogonality relation is written as:
\beq\label{so1-v2}
\sumjtom \wh_j h_a\mu_j(a)\mu_j(b\dl)=n(A)\delta_{a, b}.
\eeq
\ene

Recall that an abelian probability group is called {\it dualizable} if 
${\widehat p}_k(i,j)\geq 0$ for any $i,j,k$. In this case it can be proven that $\widehat A$ is also a probability group in which $\mu_j\dl$ is defined by $\muj\dl(a)=\mu_{j\dl}(a)=\ovl{\muj(a)}$, see \cite{hdk}.
\bp 
Let $A$ be an abelian probability group. With the above notations it follows that 
\beq\label{hatpk}
\widehat{p}_k(\jo,\jtw)=\frac{1}{n_k}\bigg(\sum_{a\in A}h_a\mujo(a)\mujtw( a)\mu_k(\adz)\bigg)
\eeq
for all $0\leq \jo,\jtw,k\leq m$.
\ep
\bpf
Evaluating both sides of Equation \eqref{hatp} at $F_k$ 
 and expanding the formula for $F_k$ from Equation \eqref{fj1} one obtains:
\begin{eqnarray*}
\widehat{p}_k(\jo,\jtw)& = & \big(\sum_l\widehat{p}_l(\jo,\jtw)\mu_l\big)(F_k) = [\mujo\star \mujtw](F_k)[\mujo\star\mujtw](F_k)\\ & =& \frac{1}{n_k}\big(\sum_{a\in A}h_a[\mujo\star\mujtw]( a)\mu_k(\adz)\big)\\ &=& \frac{1}{n_k}\big(\sum_{a\in A}h_a\mujo(a)\mujtw( a)\mu_k(\adz)\big)
\end{eqnarray*}
which gives the formula from Equation \eqref{hatpk}.
\epf
\bc 
Let $\wdht A$ be a finite probability group. With the above notations:
\beq\label{hjt}
\widehat{p}_0(\jo,\jtw)=\frac{n_{j_1}}{n(A)}\delta_{j_1, j_2^{*}}.
\eeq
\ec
\bpf For $k=0$, in Equation \eqref{hatpk} one has 
\begin{eqnarray*}
\wps_0(j_1, j_2)&=&\frac{1}{n_0}(\sum_{a\in A}h_a\mujo(a)\mujtw( a)\mu_0(\adz))\\ &= &\frac{1}{n(A)}\bigg(\sum_{a\in A}h_a\mujo(a)\mujtw( a)\bigg)\\ &=& \frac{1}{n(A)}\bigg(\sum_{a\in A}h_a\mujo(a)\mu_{\jtw\dl}( a\dl)\bigg)\\ &=& \frac{1}{n(A)}\mujo\bigg(\sum_{a\in A}h_a\mu_{\jtw\dl}( a\dl)a\bigg)\\ &\numeq{\ref{fj1}}& \frac{1}{n(A)}\mujo(n_{\jtw\dl}F_{\jtw\dl})\\ &=&\frac{n_{\jtw\dl}}{n(A)}\delta_{\jo, \jtw\dl}\\ & = & \frac{n_{\jo}}{n(A)}\delta_{\jo, \jtw\dl}.
\end{eqnarray*}
\epf
\subsection{Quotient probability group} Let $S$ be a subring of a probability ring $A$. One can define a quotient ring by the following equivalence relation:
$a\sim_Sb$ if and only if there are $s_1,s_2\in S, x\in A$ such that $m_A(1, as_1x^*)>0$ and $m_A(1, xb^*s_2)>0.
$

Define $A//S$ as the set of  equivalence classes of $\sim_S$. For an element $a \in A$ one can see that $[a]_s=[b]_S$ if and only if $u_Sau_S=u_Sbu_S$. Therefore there is a set bijection
$$
\comp[A//S] \xra{\phi} u_S\cofa u_S,\; [a]_S\sent u_Sau_S.
$$
Then $A//S$ becomes a probability group with the multiplication inherited from $u_Sau_S$ via the above isomorphism. We denote by $\bar{p}_{[c]}([a], \;[b])$ the probability structure of $A//S$. Therefore
$$[a][b]=\sum_{[c]\in A//S}\overline{p}_{[c]}([a], \;[b])[c].$$
We write shortly $[a]$ instead of $[a]_S$ when no confusion arises.
Following \cite{hdk} in the case $A$ is abelian one can show that 
\beq\label{coset}\overline{p}_{[c]}([a], \;[b])=\sum_{v \in [c]}p_v(u, s).
\eeq
with $u\in [a]_S$ and $s\in [b]_S$ arbitrarily chosen.


It was proven in \cite[Proposition 2.11]{hdk} that if $A$ is an abelian  dualizable probability group then
\beq\label{sprisom}
S^\perp \xra{\al} \widehat{A//S}, \ch \sent \al(\ch),\;\text{with}\;\al(\ch)([a]_S)=\ch(a).
\eeq
is an isomorphism of probability groups.
\subsection{}
Let $A$ be a ring which is free as $\mathbb Z$-module.
As in \cite[Chapter 3]{egno15} we define a {\it $\rgo$-ring} as a ring $A$ which is free as $\mathbb Z$-module with a special basis $B=\{b_i\}_{i \in I}$ such that 
$
b_ib_j=\sum_{k\in I}N^k_{ij}b_k
$
with $N^k_{ij}\in \rgo$. Moreover, {\it a unital $\rgo$-ring} is a $\rgo$-ring 
whose unit $1\in B$. For any unital $\rgo$-ring $(R,B)$ define a functional $\tau: R\ra \comp$ which on the basis $B$ is given by  $\tau(b):=\delta_{b,1}$.

Then a {\it $\rgo$-based ring} is a unital $\rgo$-ring togheter with an involution $i\sent i^*$ on $I$ such that the map
$$
\sumini \al_ib_i\sent \sumini \al_ib_{i^*}
$$
is a multiplicative anti-involution on $A$ and
\beq\label{tau}
\tau(b_ib_j)=\delta_{i, j^*}
\eeq
 for all $i,j\in I$.
\vsk
For the basis $B$ of a $\rgo$-based ring $(R,B)$ we can define a bilinear form $m:R\times R\ra \comp$  given by $m(\sumini \al_ib_i, \sumini \beta_ib_i)=\sumini \al_i\beta_i$.
\vsk
By \cite[Proposition 3.1.6]{egno15} the constants $N^{k^*}_{ij}$ are invariant under cyclic permutations of $\{i,j,k\}$ since $\tau$ is a trace on $R$, i.e $\tau(xy)=\tau(yx)$ for all $x,y\in R$. Therefore one can write that
$$N^{k^*}_{ij}=\tau(b_ib_jb_k)=\tau(b_jb_kb_i)=N^{i^*}_{jk}=\tau(b_kb_ib_j)=N^{j^*}_{ki}
$$
for all $i,j,k\in I$. Therefore any $\rgo$-based ring $(R,B)$ has the property that
\beq\label{br}
m(z, xy)=m(x^*, yz^*)=m(y, z^*x).
\eeq
for all $x,y, z\in R$.
\noindent
Moreover since $\tau(x)=\tau(x^*)$ for all $x\in R$ it follows that 
$$m(x, yz)=m(x^*, z^*y^*).$$
\vsk
Recall that by \cite[Proposition 3.3.9]{egno15} one has $\fp(a)=\fp(a^*)$ for all $a \in R$.
\bn{example}\label{bgr}
Given a based ring $(R, B)$ one can consider the basis $B'=\{b'\}_{b \in B}$ where $b':=\frac{b}{\fp(b)}$ and one obtains that $(B',p)$ with
\beq\label{grtprob}
p(a'b'=c')=N^c_{ab}\bigg[\frac{\fp(c)}{\fp(a)\fp(b)}\bigg]
\eeq
is a a probability group. Clearly $\cofb=R\ot_{\Z}\comp$.
\vsk
Note that in this case one has $s(a')=\fp(a)\fp(a^*)$ and 
\\ $n(B'):=\sum_{a\in B} \fp(a)\fp(a^*)=\fp(R)$.
\end{example}
\section{Fusion categories and their grothendieck groups}\label{fcgg}
Recall that a fusion category is a semisimple finite tensor  category. We refer to \cite{egno15} for the basic theory of tensor categories. 

Throughout this paper $\cc$ denotes a fusion category and $\unu$ the unit object of  a $\cc$. Given a monoidal category $\cc$ one can construct a braided fusion category $\czcc$ called the {\it monoidal centre} (or Drinfeld centre) of $\cc$, see e.g., \cite[XIII.3]{Kas} for details. The objects of $\czcc$ consist of pairs $(V, \sg_V)$ of an object $V\in \cc$ and a natural isomorphism
$\sg_{V, X}: V \ot  X \ra X \ot V $ for all $X \in \cocc$, satisfying a part of the hexagon axiom. A morphism  $f:(V, \sg_V)\ra (W, \sg_W)$ in $\czcc$
 is a morphism in $\cc$  such that $(id_X\ot f) \circ\sg_{V, X}=\sg_{W, X}\circ(f\ot id_X)$ for all objects $X$ of $ \cc$. The composition of morphisms in $\cfzcc$ is defined via  the usual composition of morphisms in $\cc$. 

Let $\cc$ be a finite tensor category and $F:\czcc\ra \cc$ be the forgetful functor. It is well known that $F$ admits a right adjoint functor $R:\cc \ra \czcc$  such that  $Z :=FR:\cc \ra \cc$ is a Hopf comonad.  Following \cite[Subsection 2.6]{scalg} one also has that 
\beq
Z(V)\simeq \int_{X\in \cc}X\ot V\ot X^*.
\eeq
If $\pi_{V;X}:Z (V)\ra X\ot V\ot X^*$ are the universal dinatural transformation associated to the above end $Z (V)$ then the counit $\eps:Z\ra \id_\cc $  is given by $\eps_V:= \pi_{V ;1} $. Moreover, using Fubini's theorem for ends, the comultiplication $\delta : Z  \ra {Z }^2$ of $Z $ is also  described in terms of the dinatural transformation $\pi$, see \cite[Subsection 3.2]{scalg}. 

The object $A:=Z(\unu)$ has the structure of  a central commutative algebra in $\cz(\cc)$. It is called the adjoint algebra of $A$.
\vsk 
Its multiplication $m_A:A\ot A \ra A$ and its unit $u_A:\unu\ra A$ are uniquely determined by  by the universal property of the end $Z$ as:
\beq 
\pi_{\unu;X} \circ u_A = \cov_X,
\eeq
\beq 
\pi_{\unu;X} \circ m_A = (\id_X \otimes \ev_X \otimes \id_{X^*} ) \circ (\pi_{\unu;X} \otimes \pi_{\unu;X})
\eeq
\noindent
Moreover $\epsu:A\ra \unu$ is a morphism of algebras, see \cite{scalg}. 

Recall that a {\it pivotal structure} $j$ on a tensor category $\cc$ is a tensor isomorphism $j:\id_\cc\ra ()^{**}$.  Given a pivotal structure one can construct for any object $X$ of $\cc$ a {\it right evaluation}
$\widetilde{ev}_X:X\ot X^*\xra{j\ot id}X^{**}\ot X^*\xra{ev_{X^*}} \unu$. 
Then for any morphism $f:A\ot X\ra B\ot X$  the right partial pivotal trace of $f$  is defined as:
\beqn
\tr_{A, B}^X:A=A\ot \unu\xra{\id \ot coev_X}A\ot X \ot X^*\xra{f \ot id} B\ot X\ot X^*\xra{\\id \ot \widetilde{ev}_X} B
\eeqn
Then the usual {\it right pivotal trace} of an endomorphism $f :X\ra X$ is obtained as a particular case for $A=B=\unu$. In particular,  the {\it right  dimension of $X$} with respect to $j$ is defined as the right trace of the identity of $X$. A  pivotal structure (or the underlying fusion category) is called {\it spherical} if
$\dim(X) = \dim(X^*)$ for all objects $X$ of $\cc$, see \cite[Section 2.2]{eno-annals}.

Given an object $X$ of $\cc$ the internal character $\mtr{ch}(X)$ is defined as the  partial pivotal trace \beqn
\mtr{ch}(X):=\tr^{X}_{A, \unu}(\rho_{X}):A\ra \unu.
\eeqn
where $\rho_X:A\ot X\ra X$ is the canonical action of $A$ on $X$, see \cite[Definition 3.3]{scalg}.

The space $\cfcc:=\hm_\cc(A, \unu)$ is called the {\it space of class functions} of $\cc$ and it is a $\comp$-algebra where the multiplication of two class functions $f, g\in \cfcc$ is defined via $f\star g:=f \circ Z(g) \circ \delta_{\unu}.$ Here the map $\delta: Z \ra  Z^2$ is the comultiplication structure of $Z$ recalled above. By \cite[Theorem 3.10]{scalg} for a pivotal fusion category $\cc$ one has that $\mtr{ch}(X\ot Y)=\mtr{ch}(X)\mtr{ch}(Y)$ for any two objects $X$ and $Y$ of $\cc$ and  $\mtr{Gr}_{\comp}(\cc)\ra \cfcc, [X]\sent \mtr{ch}(X)$ is an isomorphism of algebras. 

The  space $\cecc:= \hm_{\C}(\unu, A) $ is called {\it the space of central elements.} It is also a $\comp$-algebra where the multiplication on $\cecc$ is given by
 $a.b :=m \circ (a \ot b)$ for any $a, b \in \cecc$. 
There is a non-degenerate pairing 
$\langle\;,\; \rangle : \cfcc \times \cecc\ra \comp$, given by $ \langle f, a\rangle  \id_{\unu}= f \circ a,$ for all $f \in \cfcc$ and $a\in \cecc$. We also denote $f(a):=\langle f,\;a\rangle$. There is a right action of $\cecc$ on $\cfcc$ denoted by $\lh$  given by $f \lh b=f \circ m \circ (b\ot \id_{A})$ for all $f \in \cfcc$ and $b \in \cecc$.

Given a fusion category $\cc$ we let  $M_0, M_1, \dots M_m$  be a complete set of representatives for the isomorphism classes of simple objects of $\cc$. For brevity,  we denote $\ch_i:=\mtr{ch}(M_i)\in \cfcc$ and $d_i:=\ch_i(\unu)$ the categorical dimensions of the simple objects. Without loss of generality we may suppose that $M_0=\unu$ and therefore $d_0=1$. Moreover we denote by $i^*$ the unique index for which $M_i^*\simeq M_{i^*}$. 

Recall that also from Equation \eqref{clseq} one has 
\beq\label{clseqcc}
\sumjtom \frac{\dim(\cc)}{\dim(\cc^j)} F_j\ot F_j=\sumitom \ch_i\ot\ch_{i^*}.
\eeq


Let $\cc$ be a pivotal fusion category.  Recall that the global dimension of $\cc$ is defined as $\dimcc:=\sumitom d_i d_{i^*}\in \comp$. It is well-known that in this case $\dimcc\neq 0$. The cointegral $\lam$ of $\cc$ is defined as $\lam=\frac{1}{\dimcc}(\sumitom d_i\ch_i)\in \cfcc.$ Up to a scalar $\lam$ is the unique element of $\cfcc$ with the property that $\ch\lam=\ch(\unu)\lam$ for all $\ch \in \cfcc$.  Moreover, $\lam(u_A)=1$ and $\lam^2=\lam$. The {\it Fourier transform} of $\cc$  associated to $\lambda$ is the linear map
$
\mtc F_{\lambda}:\cecc\ra \cfcc\;\;\text{given by}\;\;a \mapsto \lambda \lh \mtc S(a)
$
where $\mtc S:\cecc\ra \cecc$ on $\cecc$ is the antipodal map of $\cecc$, see \cite[Definition 3.6]{scalg}. 


Let $\cc$ be a pivotal fusion category with commutative Grothendieck ring. Recall that $R:\cc \ra \czcc$ is the right adjoint of the  forgetful functor $F:\czcc \ra \cc$. Note that in this case by \cite[Theorem 6.6]{scalg} the object $R(\unu) \in \co(\czcc)$ is multiplicity-free.

A {\it conjugacy class} of $\cc$ is defined as a simple subobject of $R(\unu)$ in $\czcc$. Since the monoidal center $\czcc$ is also a fusion category we can write $R(\unu)=\bigoplus_{i=0}^m\mtc C_i$ as a direct sum of simple objects in $\czcc$. Thus  $\mtc C_{0},\dots, \mtc C_{m}$ are the conjugacy classes of $\C$.  Since the unit object $\unu_{\czcc }$ is always a subobject of $R(\unu)$, we may assume $\mtc C_{0} = \unu_{\czcc }$.

Let also $F_0, F_1, \dots ,F_m$ be the central primitive idempotents of $\cfcc$. We define $\mtr C_j:={\mtf}^{-1}(F_j)\in \cecc$ to be the {\it conjugacy class sums} corresponding to the  conjugacy class $\mtc C_j$.


For a fusion category $\cc$ we denote by 
$\ovl{\irr(\cc)}$ the set of elements $\{\frac{\ch_i}{d_i}\}_{i=0}^m\in \cfcc$. Since $\cfcc\simeq \mtr{Gr}_{\comp}(\cc)$ is a based ring, by Example \eqref{bgr} it follows that $\oic$ is a probability group with $\comp[\oic]=\cfcc$. By \cite[Equation 4.8]{ccc} one has that
\beq\label{njfus}
n_j=\frac{1}{F_j(\blam)}=\frac{\dim(\cc)}{\dim(\cc^j)}
\eeq
where $\blam \in \cecc$ is the idempotent integral of $\cc$.


We also denote by $\clc$ the set of central elements $\{\frac{\mtr{C}_j}{\dim(\cc^j)}\}_{j=0}^m\in \cecc$.  We denote by $\widehta{\cfcc}$ the dual complex algebra $\comp[\widehta{\oic}]$.
\vsk
By the proof of  \cite[Theorem 6.12]{scalg} for all $\ch \in \cfcc$ one has that 
\beq\label{sh}
\muj(\ch)=\ch(\frac{\barcj}{\dim(\cc^j)}).
\eeq


Since $\oic$ is an abelian probability group orthogonality relations form Equations \eqref{so1} and \eqref{so1-v2} can be written as
\beq\label{orthof}
\sum_{ k=0}^m\frac{\ch_i(\barck)\ch_{j^*}(\barck)}{\dim(\cc^k)}=\delta_{i,j}\dimcc
\eeq
and respectively 
\beq\label{orthos}
\sum_{ i=0}^m\ch_i(\barcl)\ch_{i^*}(\barck)=\delta_{l,k}{\dim(\cc^k)}\dimcc
\eeq
\bc\label{div}
Let $\cc$ be a fusion category with a commutative Grothendieck ring.
$$
\frac{\dim(\cc)}{\dim(\cc^j)} \in \mathbb A
$$
\ec
\bpf Since as above $\muj(\ch_i)=\frac{\ch_i(\barcj)}{\dim(\cc^j)}$ (see the proof of \cite[Theorem 6.12]{scalg}) the first orthogonality relation \eqref{orthof} can be written as
$\sum_{i=0}^m\mu_k(\ch_i)\mu_k(\ch_{i^*})=\frac{\dimcc}{{\dim(\cc^k)}}$.  On the other hand it is known that $\mul(\ch_i)$ are cyclotomic integers for all  $0\leq l \leq m$.
\epf
\bl For any $\ch \in \cfcc$ and $z, z'\in \cecc$ one has
\beq\label{zz}
d(\ch)\ch(zz')=\ch(z)\ch(z')
\eeq
\el
\bpf If $z=\sumitom z_iE_i$ and $z'=\sumitom z'_iE_i$ then $zz'=\sumitom z_iz'_iE_i$
Then $\ch_i(z)\ch_i(z')=z_iz'_id_i^2$ and $d(\ch_i)\ch_i(zz')=d_i^2z_iz'_i=\ch_i(z)\ch_i(z')$.
\epf

Next theorem extends the result obtained in \cite{zz} for semisimple Hopf algebras.
\bt\label{cec}
Let $\cc$ be a fusion category. Then $\clc$ is an abelian {generalized hypergroup} and $\clc\simeq\wht{\oic}$ as generalized hypergroups.
More precisely, the map
\beq\label{dzc}
\wht{\oic} \xra{\simeq} \clc, \;\muj\sent \frac{\mtr C_j}{\dim(\cc^j)}
\eeq
is bijective and induces an isomorphism of $\comp$-algebras
$\wht{\cfcc} \xra{\simeq} \cecc.$
\et

\bpf
Clearly this map is bijective, it remains to show that it induces  an isomorphism of $\kk$-algebras $\al: \wht{\cfcc} \xra{} \cecc,\; \muj\sent \frac{\mtr C_j}{\dim(\cc^j)}$.

One has to prove that $\al(\mui\muj)=\al(\mui)\al(\muj)$ for all $i,j$. Suppose as above that $\mui \muj=\sum_k{\wdht p}_k(i,j)\muk.$ Evaluating at $\frac{\ch_s}{d_s}$ this identity it follows that
$$
[\mui \muj](\frac{\ch_s}{d_s})=\sum_k{\wdht p}_k(i,j)\muk(\frac{\ch_s}{d_s})=\frac{\ch_s}{d_s}(\sum_k{\wdht p}_k(i,j)\frac{\barck}{\dim(\cc^k)})
$$
On the other hand 
\begin{eqnarray*}
[\mui \muj](\frac{\ch_s}{d_s})&=&\mui(\frac{\ch_s}{d_s})\muj(\frac{\ch_s}{d_s})=\frac{\ch_s}{d_s}(\frac{\barci}{\dim(\cc^i)})\frac{\ch_s}{d_s}(\frac{\barcj}{\dim(\cc^j)})
\\ &\numeq{\ref{zz}}& \frac{\ch_s}{d_s}(\frac{\barci}{\dim(\cc^i)}\frac{\barcj}{\dim(\cc^j)})
\end{eqnarray*}
Thus 
$$
\sum_k{\wdht p}_k(i,j)\frac{\barck}{\dim(\cc^k)}=\frac{\barci}{\dim(\cc^i)}\frac{\barcj}{\dim(\cc^j)}
$$
which shows that $\al(\mu_i\muj)=\al(\mu_i)\al(\muj)$
\epf
\subsection{Proof of Theorem \ref{burnside}}\label{proof}
In this section we give an analogue of Burnside formula for the structure constants. Since $\{\barcj\}_j$ form a $\comp$-linear basis for $\cecc$ one has that 
\beq\label{defstrc}
\barcjo\barcjtw=\sum_{l=0}^m c^l_{\jo,\jtw}\barcl
\eeq
form some scalars $c^l_{\jo,\jtw}\in \comp$. These scalars are called {\it the structure constants of $\cc$.} The last equation from the proof above shows that
$$
c^k_{ij}=\frac{\dim(\cc^i)\dim(\cc^j)}{\dim(\cc^k)}{\wdht p}_k(i,j).
$$
Combined with Equation \eqref{hatpk} one obtains the result of  Theorem \ref{burnside}.

\section{Proof of Theorem \ref{main2}}\label{braided}
A braided category is called {\it premodular} if it  has a spherical structure. Equivalently, this  is a ribbon fusion category, that is, a fusion category equipped with a braiding and a twist
(also called a balanced structure), see \cite{eno-annals}.

Then $\cc$ is called {\it pseudo-unitary} if $\fpcc =\dimcc$. If such is the case,
then by \cite[Proposition 8.23]{eno-annals}, $\cc$ admits a unique spherical structure with respect to which the categorical dimensions of simple objects are all positive. It is called the canonical spherical structure. For this structure, the categorical dimension of an object coincides with its Frobenius-Perron dimension, i.e. $\fp(X)=\dim(X)$ for any object $X$ of $\cc$. 

Recall \cite{eno-adv} that a fusion category $\cc$ is called {\it weakly integral} if its Frobenius-Perron dimension is an integer. If $\cc$ is such a  fusion category then $\cc$ is pseudo-unitary by \cite[Proposition 8.24]{eno-annals}. 

Moreover, if $\cc$ is weakly integral then by \cite[Proposition 8.27]{eno-annals} the dimensions of  simple objects in $\cc_{ad}$ are integers. Since the conjugacy classes $\cc^j$ are sum of simple object of the adjoint subcategory it follows that $\dimccj$ are integers.

Recall that premodular fusion category $\cc$ is called {\it modular} if its $S$-matrix is non-degenerate.
\bt\label{dmc}
Let $\cc$ be a modular fusion category. Then $\oic$ is a dualizable probability group and
$$\oic\simeq \widehat{\oic}$$
as probability groups.
\et
\bpf
By \cite[Example 6.14]{scalg}  the map $f_Q:\cfcc \ra \cecc, \ch_i\sent \sumjtom \frac{\sij}{d_j}e_j$ is an isomorphism of algebras. Moreover from Proposition \cite[Theorem 6.5]{ccc} one has that $f_Q(\frac{\ch_i}{d_i})=g_i=\frac{\mtr C_i}{\dim(\cc^i)}$ which shows that $f_Q$ sends bijectively $\oic$ in $\clc$. This shows that $\clc$ is also a probability group and $f_Q:\oic \ra \clc$ is an isomorphism of probability groups. On the other hand by Theorem \ref{cec} one has that in this situation $\clc\simeq \widehat{\oic}$ as probability groups.
\epf
Let $\cc$ be a pivotal fusion category with a commutative Grothendieck ring $\cfcc\simeq \mtr{Gr}_{\kk}(\cc)$. Let also $F_0, F_1, \dots ,F_m$ be the primitive central idempotents of $\cfcc$ and  $\mu_0, \mu_1, \dots \mu_m$ be their corresponding characters on $\cfcc$. Therefore
$$\mui:\cfcc\ra \comp,\;\mui(F_j)=\delta_{i,j}.$$
We also denote by $\cc^0, \cc^1,\dots \cc^m$ the conjugacy classes of $\cc$ corresponding in this order to the primitive idempotents $F_0, F_1, \dots F_m$.

Let also $M_0, M_1, \dots M_m$ be a complete set of representatives for the isomorphism classes of simple objects of $\cc$. As above, without loss of generality we may assume that $M_0=\unu$ is the unit object of $\cc$.

As above we denote by $V_0, V_1,\dots V_r$  a complete set of  simple objects of $\czcc$ and by $\we_0, \we_1, \dots \we_r$ their associated primitive idempotents in $\mtr{CE}(\czcc)$. We denote also by $\wch_0, \wch_1, \dots \wch_r$ the characters associated to $V_0, V_1, \dots V_r$ and let $\wzd_0,\;\wzd_1,\;\dots \wzd_r$ be their quantum dimensions. Therefore $\wzd_s=\wch_s(\unu)$ for all $s$.
We may also assume that $V_0=\unu$ is the unit object of $\czcc$. Without loss of generality we may also suppose that $V_i=\cc^i$ for any $0\leq i \leq m$. Since the Drinfeld map $F_Q:\cf(\czcc)\ra \mtr{CE}(\czcc)$ is bijective it follows that $\wf_j:=F_Q^{-1}(\we_j)$ is a complete set of primitive orthogonal idempotents of $\cf(\czcc)$.

Let also $F:\czcc\ra \cc$ be the forgetful functor. It is well known, see \cite[Lemma 8.49]{eno-annals}, that the induced map $\res=F_*:\cf(\czcc)\ra \cfcc$ is surjective.

Moreover, Ostrik showed in \cite[Theorem 2.13]{O3} that for any primitive idempotent $F_j\in \cfcc$ of a spherical category $\cc$ there is a unique primitive idempotent $\wtf_{\sgj}\in \cf(\czcc)$ whose restriction to $\cfcc$ is $F_j$ and moreover $V_{\sgj}=\cc_{\sg(j)}$ is a conjugacy class of $\cc$ with 
$$
\dim(V_{\sgj}) =\dim(\ccj),\;\text{i.e}\; \wzd_{\sg(j)}=\wzd_j.
$$
Thus 
$F_*(\wf_{\sg(j)})=F_j,\;\tetx{for all}\; 0\leq j \leq m.$ Note also that $F_*(\wf_s)=0$ for $s\neq \sg(j)$ for some $j$. We denote by $\sz_{ij}$ the $S$-matrix of the braided category $\czcc$. 
\subsection{} For the rest of this section we suppose that $\cc$ is a braided spherical fusion category, i.e a premodular category. Since in this case $\czcc$ is a modular tensor category  it follows by \cite[Theorem 3.1.12]{BaKi} that the irreducible characters of $\cf(\czcc)$ are indexed by the simple object of $\czcc$. More precisely, if $V_i$ is a simple object of $\czcc$ then
$$
\wmu_i:\cf(\czcc)\ra \comp,\; \wmui([V_j])=\frac{\sz_{ij}}{d_i}
$$
is an algebra homomorphism.

%

By \cite[Section 2.10]{dgno2} there is also a braided tensor functor 
\beq\label{iota}
\iota:\cc \hookrightarrow \czcc, X\sent (X, c_{X,-}).
\eeq
that is fully faithful.
It follows that $\iota(M_i)\simeq V_{\tilde{i}}$ for  some $0\leq \tilde{i}\leq r$.  Note that $\{\tilde 0, \tilde 1, \dots ,\tilde m\}\cap \{0,1,\dots,  m\}=\{0\}$. Indeed, for $i>1$ $\iota(M_i)$ cannot be a conjugacy class since $M_i=F(\iota(M_i))$ does not contain the unit object $\unu$ of $\cc$. Since $F\circ \iota=\id_\cc$ note that 
$F_*(\wch_{\tilde t})=\ch_t,\;\text{for all}\; 0 \leq t \leq m$.

Consider the inclusion $\cfcc\subseteq \cf(\czcc)$ induced by the inclusion functor of \eqref{iota}. For any $0\leq j \leq m$ we define a class of characters
$$
\ca_j:=\{\wmui\in \widehat{\cf(\czcc)}\;|\;\wmui|_{\cfc}=\muj\}.
$$
\bp\label{sigma} Let $\cc$ be a premodular category.
With the above notations one has $\sg(j)\in \ca_j$ for any  $0\leq j \leq m$. 
\ep
\bpf 
Following \cite[Example 6.14]{scalg}, inside $\cf(\czcc)$ one can write that
$$
\wch_{\tilde i}=\sum_{l=0}^{r} \frac{\sz_{\tilde i l}}{\wzd_l} \wf_l, 
$$
for all $0\leq i\leq m$.
Applying the morphism $F_*:\cf(\czcc)\ra \cfcc$ induced by the forgetful functor $F$ one obtains that
$$
\ch_i=\sum_{j=0}^{r} \frac{\sz_{\tilde i j}}{\wzd_{j}}F_*(\wf_j)=\sum_{j=0}^m\frac{\sz_{\tilde i \sg(j)}}{\widehat{d}_{\sg(j)}}F_j
.$$
This implies that $\muj(\ch_t)=\frac{\sz_{s \sg(j)}}{d_{\sg(j)}}=\wmu_{\sgj}(\wch_{\tilde t}),$ thus $\wmu_{\sgj}\in \ca_j$.
\epf
\subsection{Proof of the main Theorem \ref{main2}}
We denote by $P_v(s, u)$ the probability structure of $\ovl{\irr(\czcc)}$ and by $\wpm_k(i,j)$ the probability structure of the dual abelian probability group $\widehta{\ovl{\irr(\czcc)}}$.	By Example \eqref{bgr} we have 
$$\wpm_v(s,u)=P_v(s, u)=\frac{\widehta{N}^v_{s, u}\widehat{d}_v}{\widehat{d}_s\widehat{d}_u}$$
\bpf
Since $\cc$ is braided we have an inclusion of braided categories
$$
\cc \hookrightarrow \czcc, X\sent (X, c_{X,-}).
$$
This shows that $\cfcc \subseteq \cf(\czcc)$ and therefore $\oic \leq \overline{\irr({\czcc)}}$ is a probability subgroup.  By \cite[Proposition  2.11 ]{hdk} it follows that $\oic$ is a dualizable probability group and
\beq\label{isom}
\comp[\widehat{\overline{\irr({\cc)}}}]=\widehat{\cfcc}\simeq \widehat{\cf(\czcc)}//\cfcc^{\perp}
\eeq
where
$
\cfcc^\perp=\{\wmu_j\in \widehat{\cf(\czcc)}\;|\; {\wmu_j}(\wch_s)=\fp(\ch_s), \;\tetx{for all}\;0\leq s \leq m\}.
$
Moreover, under the isomorphism \eqref{sprisom} one has 
$$[\wmu_u]_{\cfcc^\perp}=[\wmu_v]_{\cfcc^\perp} \iff u, v \in \ca_j,\; \text{for some}\;j.$$ 
\vsk By Equation \eqref{coset} one has
\beq\label{wprat}
\wdhat{p}_k(i,j)=\sum_{v\in \ca_k}\wdhat{P_v}(s,u)=\sum_{v\in \ca_k}P_v(s,u)=\sum_{v\in \ca_k}\frac{\wdhat{N^v_{su}}\wdhat{d_v}}{\wdhat{d_s}\wdhat{d_u}}.
\eeq
\vsk 
Since $\mtr C_i=\dim(\cc^i)g_i$ it follows that
$$
\mtr C_i\mtr C_j=\sum_{k=1}^m\bigg(\sum_{v\in \ca_k}\frac{{\wdht N}^v_{su}\wdhat{d_v}\dim(\cc^i)\dim(\cc^j)}{\wdhat{d_s}\wdhat{d_u}\dim(\cc^k)}\bigg)\mtr C_k
$$
By Proposition \ref{sigma} one has $\sg(i)\in \ca_i$. Therefore if $s=\sg(i)$ and $u=\sg(j)$, then
$$
\mtr C_i\mtr C_j=\sum_{k=1}^m\bigg(\sum_{v\in \ca_k}\frac{{\wdht N}^v_{\sg(i)\sg(j)}\widehat{d}_v}{\dim(\cc^k)}\bigg)\mtr C_k
$$
since $\wdhat{d}_{\sg(i)}=\dim(\cc^i)$ and $\wdhat{d}_{\sg(j)}=\dim(\cc^j)$.
\vsk By Corollary \ref{div} one has $\dim(\cc^k)\big|\dim(\cc)$ and therefore 
$$
\dim(\cc)c^k_{ij}=\frac{\dim(\cc)}{\dim(\cc^k)}\bigg (\sum_{v\in \ca_k}{\wht N}^v_{\sg(i)\sg(j)}\widehat{d}_v \bigg)
$$
is an algebraic integer. Moreover, if $\cc$ is weakly integral then by the above formula, the same numbers are  non-negative integers. Indeed, by \cite[Theorem 3.10]{NG} one has  $\widehat{d}_v$ is integer for any $v \in \ca_k$ since $\cc_{\sgi}$ and $\cc_{\sg(j)}$ have integer dimensions. Moreover $\frac{\dim(\cc)}{\dim(\cc^k)}$ is a rational number which is also  an algebraic integer, therefore an integer.
\epf
\br
Taking $j=i^*$ in the proof of the above theorem and $u=s^*$ it follows that
$
\frac{1}{\dim(\cc^i)}=\wdhat{p}_1(i,i^*)=\sum_{v\in \ca_1}\frac{{\widehat{N}^v}_{ss'}\widehat{d}_v}{{\wdhat{d}_s}^2}
$
i.e. $\dim(\cc^i)\big|\widehat{d}_s^{\;2}$ for all $s \in \ca_i$.

In particular, if $\wdhat{d}_s=1$ and $s \in \ca_j$ it follows that $\dim(\cc^j)=1$. 
\er
\bibliographystyle{alpha}
\bibliography{scbf-8}
\ed

\bibliographystyle{alpha}
\bibliography{scbf8}
\ed